\documentclass{amsart}
\usepackage{amsfonts, amsmath, amssymb, amsthm, amsxtra, latexsym}

\theoremstyle{plain}

\newtheorem{theorem}{Theorem}[section]
\newtheorem{proposition}[theorem]{Proposition}
\newtheorem{lemma}[theorem]{Lemma}

\theoremstyle{definition}

\newtheorem{definition}[theorem]{Definition}

\numberwithin{equation}{section}
\newcommand {\C}{{{\mathbb C}}}
\newcommand {\N}{{{\mathbb N}}}
\newcommand {\R}{{{\mathbb R}}}

\newcommand {\Nn}{{\N^n}}
\newcommand {\Rn}{{\R^n}}

\newcommand {\measures}{{{\mathcal M}^*}}
\newcommand {\measuresC}{{{\mathcal M}^*_C}}
\newcommand {\trigspan}{ { \Span_\C \, \{e_{i\lambda}\mid \lambda\in S\} } }
\newcommand {\wtilde}{{\widetilde{w}}}

\def\Span{\mathop{\rm Span}\nolimits}

\begin{document}
%%%%%%%%%%%%%%%%%%%%%%%%%%%%%%%%%%%%%%%%%%%%%%%%%%%%%%%%%%%%%%%%%%%%%%%%%%%%%%%%%%%%%%%

%%%%%%%%%%%%%%   Top matter %%%%%%%%%%%%%%%%%%%%%%%%%%%%%%%%%%%%%%%%%%%%%%%%%%%%%%%%%%%

\title[Determinate measures]{Determinate multidimensional measures, \\
       the extended Carleman theorem\\ and quasi-analytic weights}
             % Linebreaks by \\, leave blank if not needed

\translator{}
             % Linebreaks by \\, leave blank if not needed

\dedicatory{}
            % Linebreaks by \\, leave blank if not needed

\author[Marcel~de~Jeu]{Marcel~de~Jeu}
             % Comment out if not needed

\begin{abstract}

We prove in a direct fashion that a multidimensional probability measure $\mu$ is determinate if the higher dimensional analogue of Carleman's condition is
satisfied. In that case, the polynomials, as well as certain proper subspaces of the
trigonometric functions, are dense in all associated $L_p$-spaces for $1\leq p<\infty$. In
particular these three statements hold if the reciprocal of a quasi-analytic weight has finite
integral under $\mu$. We give practical examples of such weights, based on their classification.

As in the one dimensional case, the results on determinacy of measures supported on $\Rn$ lead
to sufficient conditions for determinacy of measures supported in a positive convex cone, i.e.\ the higher dimensional analogue of determinacy in the sense of Stieltjes.

\end{abstract}
             % Comment out if not needed

\date{}
             % Footnote on page 1, first item
             % Leave blank if not needed

             % Footnote on page 1, second item
             % Enter as {Primary primarynumber; Secondary secundarynumbers}
             % Leave blank if not needed

             % Footnote on page 1, third item
             % Enter the keywords, separated by commas, no period at the end.
             % Leave blank if not needed

\thanks{During the preparation of this paper the author was partially supported by a PIONIER grant of the Netherlands Organisation for Scientific Research (NWO)}
             % Footnote on page 1, fourth item.
             % Do not specify linebreaks
             % Leave blank if not needed

\address{M.F.E.~de~Jeu\\
         Korteweg-de Vries Institute for Mathematics\\
         University of Amsterdam\\
         Plantage Muidergracht~24\\
         1018~TV~Amsterdam\\
         The~Netherlands}
             % Comment out if not needed

%\curraddress{}
             % Linebreaks by \\, leave blank if not needed.
             % Command is not recognized????? (October 8, 2000)

\email{mdejeu@science.uva.nl}
             % Comment out if not needed

\urladdr{}
             % Leave blank if not needed

\maketitle

%%%%%%%%%%%%%%   Start of the text %%%%%%%%%%%%%%%%%%%%%%%%%%%%%%%%%%%%%%%%%%%%%%%%%%%%
\section{Introduction and overview}

We will be concerned with determinacy and density results for probability measures on $\Rn$
for a fixed $n$. Establishing notation, let $(\,.\,)$ be the standard inner
product on $\Rn$ with corresponding norm $\Vert\,.\,\Vert$. For $\lambda\in\Rn$
define $e_{i\lambda}:\Rn\mapsto\C$ by $e_{i\lambda}(x)=\exp i(\lambda,x)$
($x\in\Rn$). We let $\measures$ be the set of all positive Borel measures $\mu$
on $\Rn$ such that
\begin{equation*}
\int_\Rn \Vert x\Vert^d\,d\mu(x)<\infty
\end{equation*}
for all $d\geq 0$. A measure $\mu$ is said to be determinate if $\mu\in\measures$ and if $\mu$ is uniquely determined in $\measures$ by the set
of integrals
\begin{equation*}
\int_\Rn P(x)\,d\mu(x)
\end{equation*}
of all polynomials $P$ on $\Rn$.

There are several results in the literature concerning the determinacy of elements of
$\measures$ and the related matter of the density of the polynomials in the associated
$L_p$-spaces. Connections are furthermore known between these properties for a multidimensional
measure and the corresponding properties for its one dimensional marginal distributions. We
refer to \cite{Berg1,Berg2} for an overview of the field.

These results in the literature yield sufficient conditions for a measure to be
determinate. The resulting criteria are however not always easy to apply, since
they tend to ultimately involve the computation of moment sequences. Given a
particular measure such a computation need not be an attractive task.

In this paper on the contrary we establish an \emph{integral} criterion of some generality to
conclude that a measure is determinate. A criterion of this type is evidently easier to apply.
Moreover, if this criterion is satisfied, then the polynomials are dense in the associated
$L_p$-spaces for finite $p$, and the same holds for $\trigspan$ for any subset $S$ of $\Rn$
which is somewhere dense, i.e.\ is such that its closure $\overline S$ has non-empty interior.
Our criterion is established along the following lines.

We first prove that a multidimensional probability measure is determinate and that the density
results hold as described above if the higher dimensional analogue of Carleman's condition is
satisfied. This should, analogously to \cite{BergChristensen2}, be compared with the classical
one dimensional Carleman theorem, which asserts determinacy but is not concerned with density.
We will in arbitrary dimension refer to the total conclusion of the determinacy and the density
as described above as the \emph{extended} Carleman theorem.

Our proof of the extended Carleman theorem is based on a result on multidimensional
quasi-analytic classes. It is a ``direct'' proof and close in spirit to the classical proof of
the one dimensional Carleman theorem as in e.g.\ \cite{Koosis}. We will also indicate an
alternative derivation, based on the recent literature. This alternative derivation however is
considerably less direct than our approach.

Having established the extended Carleman theorem, we subsequently note that a measure satisfies
the necessary hypotheses if the reciprocal of a so called quasi-analytic weight has finite
integral. Such (multidimensional) weights are defined and studied systematically in
\cite{deJeu}. The sufficiency of the aforementioned integral condition on the measure is then in
fact almost trivial, given the definition of these weights in terms of divergent series as in
Section \ref{sec:qaw}. Verifying this divergence for a particular weight is however in general
not an easy computation, but---and this is the crucial point---quasi-analytic weights can
alternatively be characterized by the divergence of certain integrals, which is on the contrary
usually a rather straightforward condition to verify. Using these two equivalent
characterizations, we are thus finally led to results in the vein of the following theorem.

\begin{theorem}\label{thm:introthm}
Suppose $R>0$ and a non-decreasing function $\rho:(R,\infty)\mapsto\R_{\geq 0}$ of class $C^1$
are such that
\begin{equation*}
\int_{R}^\infty\frac{\rho(s)}{s^2}\,ds=\infty.
\end{equation*}
If $\mu$ is a positive Borel measure on $\Rn$ such that
\begin{equation*}
\int_{\Vert x\Vert>R} \exp{\left(\int_{R}^{\Vert x\Vert}\frac{\rho(s)}{s}\,ds\right)}\,d\mu(x)<\infty,
\end{equation*}
then $\mu$ is determinate. Furthermore, the polynomials and $\trigspan$ are then
dense in $L_p(\Rn,\mu)$ for all $1\leq p<\infty$ and for every subset $S$ of
$\Rn$ which is somewhere dense.
\end{theorem}

A particular case is obtained by choosing $\rho(s)=\epsilon s$ for some
$\epsilon>0$. Then one sees that
\begin{equation*}
\int_\Rn \exp{\left(\epsilon \Vert x\Vert\right)}\,d\mu(x)<\infty
\end{equation*}
is a sufficient condition; this is a classical type of result. However, a measure
$\mu\in\measures$ is now also seen to be determinate, and the polynomials and spaces $\trigspan$
for somewhere dense $S$ are dense in the associated $L_p$-spaces for finite $p$, if e.g.\
\begin{equation}\label{eq:logintfinite}
\int_{a_2 \Vert x\Vert>2} \exp{\left(\frac{a_1 \Vert x\Vert}{\log a_2\Vert
x\Vert}\right)}\,d\mu(x)<\infty
\end{equation}
for some $a_1,\,a_2>0$. This is a substantially weaker condition. In Section
\ref{sec:Hamburger} we will give some additional and even more lenient
sufficient conditions, formulated in terms of elementary functions as above. It
will also become apparent that the integrand need not be radial as in Theorem
\ref{thm:introthm}. Although such radial integrands may be sufficient for most
applications, this is not the most general situation in which our results apply.
We return to the possible consequences of this observation in Section~\ref{sec:closing}.

In the discussion so far we considered what can be called determinacy in the sense of Hamburger
in arbitrary dimension, i.e.\ the question whether a measure on $\Rn$ is determined by its
integrals of the polynomials, without any restriction on its support. Naturally, in the one
dimensional case the question of determinacy has also been studied under the condition that the
support of the measure is contained in the interval $[0,\infty)$. This determinacy in the sense
of Stieltjes has an analogue in arbitrary dimension, by asking whether a measure on $\Rn$ is
determined by its integrals of the polynomials, under the assumption that its support is
contained in a given positive convex cone with the origin as vertex. The simultaneous
distribution of non-negative random variables provides an obvious practical example. To
facilitate the formulation, we adapt the following terminology.

\begin{definition}\label{def:cdet}
Let $\{v_1,\ldots,v_n\}$ be a basis of $\Rn$ and let $C=\R_{\geq 0}\cdot v_1+\ldots\R_{\geq
0}\cdot v_n$ be the corresponding positive convex cone. Let $\mu\in\measures$ be supported in
$C$. Then $\mu$ is \emph{$C$-determinate} if a measure $\nu\in\measures$, which is also supported in $C$ and is such that
\begin{equation*}
\int_C P(x)\,d\nu(x)=\int_C P(x)\,d\mu(x)
\end{equation*}
for all polynomials $P$ on $\Rn$, is necessarily equal to $\mu$.
\end{definition}

As in the one dimensional case, sufficient conditions for determinacy in the sense of Hamburger
imply sufficient conditions for $C$-determinacy. The density results do not transfer in general.
Concentrating on radial integrands again, we thus obtain the following result.

\begin{theorem}\label{thm:introthmStieltjes}
Let $\{v_1,\ldots,v_n\}$ be a basis of $\Rn$ and let $C=\R_{\geq 0}\cdot v_1+\ldots\R_{\geq
0}\cdot v_n$ be the corresponding positive convex cone. Let $\mu$ be a positive Borel measure on
$\Rn$ which is supported in $C$.

Suppose $R>0$ and a non-decreasing function $\rho:(R,\infty)\mapsto\R_{\geq 0}$ of class $C^1$
are such that
\begin{equation*}
\int_{R}^\infty\frac{\rho(s)}{s^2}\,ds=\infty
\end{equation*}
and
\begin{equation*}
\int_{\sqrt{\Vert x\Vert}>R} \exp{\left(\int_{R}^{\sqrt{\Vert
x\Vert}}\frac{\rho(s)}{s}\,ds\right)}\,d\mu(x)<\infty.
\end{equation*}
Then $\mu$ is $C$-determinate.
\end{theorem}

Aside, we mention that under an additional condition one can conclude that $\mu$ is actually
determinate, as will be discussed in Section~\ref{sec:Stieltjes}.

As a consequence of the theorem, if $\mu\in\measures$ is supported in $C$, and if e.g.\
\begin{equation*}
\int_C \exp{\left(\epsilon \sqrt{\Vert x\Vert}\right)}\,d\mu(x)<\infty
\end{equation*}
for some $\epsilon>0$, or if
\begin{equation}\label{eq:sqrtlogintfinite}
\int_{a_2 \sqrt{\Vert x\Vert}>2} \exp{\left(\frac{a_1 \sqrt{\Vert x\Vert}}{\log a_2\Vert
x\Vert}\right)}\,d\mu(x)<\infty
\end{equation}
for some $a_1$, $a_2>0$, then $\mu$ is $C$-determinate.

To conclude this introductory discussion, we first of all mention that for quite a few (one
dimensional) common distributions occurring in practice the determinacy or non-determinacy is
known; see e.g. \cite{PakesHungWu} for a number of examples. It seems that many of the known
positive results on determinacy in one dimension follow from condition \eqref{eq:logintfinite}
in the Hamburger case or \eqref{eq:sqrtlogintfinite} in the Stieltjes case, the latter possibly
combined with the result as discussed at the end of Section~\ref{sec:Stieltjes}.

Secondly, let us note that the typical practical sufficient condition from which we conclude
determinacy in this paper is the integrability of a function of a suitable type. The
``underlying'' reason for this determinacy is a Carleman-type criterion, which is satisfied as a
consequence of this integrability. It is an interesting problem to determine a set of functions
with the property that a measure satisfies such a Carleman-type criterion \emph{precisely} if a
function in this set is integrable. These matters are addressed in \cite{Hoffmann} for the one
dimensional case.

This paper is organized as follows.

In Section~\ref{sec:Carleman} we establish the extended Carleman theorem.

Section~\ref{sec:qaw} is a preparation for Sections~\ref{sec:Hamburger} and \ref{sec:Stieltjes}.
It contains the definition of quasi-analytic weights, their classification and main properties,
referring to \cite{deJeu} for proofs.

In Section~\ref{sec:Hamburger} the results of Sections~\ref{sec:Carleman} and \ref{sec:qaw} are
put together, resulting in integral criteria for determinacy (without restrictions on the
support) and the density results.

Section~\ref{sec:Stieltjes} is concerned with determinacy in the sense of Stieltjes, i.e.\ with
$C$-determinacy. Integral criteria are obtained and a condition is discussed under which one can
conclude determinacy, rather than just $C$-determinacy.

Section~\ref{sec:closing} contains a tentative remark on the possibility of the
existence of distinguished marginal distributions.

\section*{Acknowledgments.}
It is a pleasure to thank Christian Berg for a helpful exposition on the subject and useful
comments on a previous version of the paper.

\section{The extended Carleman theorem}\label{sec:Carleman}

In this section we establish the extended Carleman theorem. The determinacy of the measure, the
density of the polynomials and the density of the spaces $\trigspan$ are all seen to be closely
related, since they all ultimately rest on the following theorem on multidimensional
quasi-analytic classes.

\begin{theorem}\label{thm:DCtheorem}
For $j=1,\ldots,n$ let $\{M_j(m)\}_{m=0}^\infty$ be a sequence of non-negative
real numbers such that
\begin{equation*}
\sum_{m=1}^\infty \frac{1}{M_j(m)^{1/m}}=\infty.
\end{equation*}
Assume that $f:\R^n\mapsto\C$ is of class $C^\infty$ and that there exists
$C\geq 0$ such that
\begin{equation*}
\left\vert \frac{\partial^\alpha f}{\partial \lambda^\alpha}(\lambda)\right\vert\leq C\prod_{j=1}^n M_j(\alpha_j)
\end{equation*}
for all $\alpha\in\Nn$ and all $\lambda\in \R^n$. Then, if $
\frac{\partial^\alpha f}{\partial \lambda^\alpha}(0)=0$ for all $\alpha\in\Nn$,
$f$ is actually identically zero on $\R^n$.
\end{theorem}

A proof by induction, starting from the Denjoy-Carleman theorem in one
dimension, can be found in \cite{deJeu}. The result in [loc.cit.] is in fact
somewhat stronger than the statement above. A slightly weaker version on the
other hand, which is however not entirely sufficient in our situation, can already be
found in \cite{Hryptun}. The proof in \cite{Hryptun} is more complicated than the
proof in \cite{deJeu}, but in \cite{Hryptun} the necessity of the hypotheses is
investigated as well.

For reference purposes we state the following elementary fact, the verification
of which is omitted.

\begin{lemma}\label{lem:elementary} Let $\{a(m)\}_{m=1}^\infty$ be a non-negative
non-increasing sequence of real numbers. If $k$ and $l$
are strictly positive integers, then $\sum_{m=1}^\infty a(km)=\infty$ if and
only if $\sum_{m=1}^\infty a(lm)=\infty$.
\end{lemma}

\begin{theorem}[Extended Carleman theorem]\label{thm:Carleman}
Let $\mu\in\measures$ and suppose $\{v_1,\ldots,v_n\}$ is a basis of $\Rn$. For
$j=1,\ldots,n$ and $m=0,1,2,\ldots$ define
\begin{equation*}
s_j(m)=\int_\Rn (v_j,x)^m\,d\mu(x).
\end{equation*}
If each of the sequences $\{s_j(m)\}_{m=1}^\infty$ ($j=1,\ldots,n$) satisfies
Carleman's condition
\begin{equation}\label{eq:Carleman}
\sum_{m=1}^\infty\frac{1}{s_j(2m)^{1/2m}}=\infty,
\end{equation}
then $\mu$ is determinate. Furthermore, the polynomials and $\trigspan$ are then
dense in $L_p(\Rn,\mu)$ for all $1\leq p<\infty$ and for every subset $S$ of
$\Rn$ which is somewhere dense.
\end{theorem}

\begin{proof}{Proof} Using the obvious fact that a linear automorphism of $\Rn$ induces an automorphism of
the polynomials and a permutation of the spaces of trigonometric functions $\trigspan$ for
somewhere dense $S$, one sees easily that we may assume that $\{v_1,\ldots,v_n\}$ is the
standard basis of $\Rn$. One also verifies that we may assume in addition that $\mu(\Rn)=1$.
Under these two assumptions we turn to the proof.

Treating the determinacy of $\mu$ first, we write $\mu_1=\mu$ and we suppose
that $\mu_2\in\measures$ is a probability measure on $\Rn$ with the same
integrals of the polynomials as $\mu_1$. Let $\nu=\frac{1}{2}(\mu_1+\mu_2)$ and
introduce, in the usual multi-index notation:
\begin{equation*}
t(\alpha)=\int_\Rn \vert x\vert^\alpha \,d\nu\quad(\alpha\in\Nn).
\end{equation*}
Let
\begin{equation}\label{eq:marginalsequences}
t_j(s)=\int_\Rn \vert x_j\vert^s \,d\nu\quad(j=1,\ldots,n;\,s\geq 0).
\end{equation}
As is well known, the H\"older inequality and the fact that $\nu(\Rn)=1$ imply
that
\begin{equation}\label{eq:increasing}
t_j(s_1)^{1/s_1}\leq t_j(s_2)^{1/s_2}
\end{equation}
for $j=1,\ldots,n$ and $1\leq s_1\leq s_2<\infty$. In addition, regarding $\vert
x\vert^\alpha=\prod_{j=1}^n \vert x_j\vert^{\alpha_j}$ as a product of $n$ elements of
$L_n(\Rn,\nu)$, the generalized H\"older inequality \cite[VI.11.1]{DunfordSchwartz} yields
\begin{equation}\label{eq:Holder}
t(\alpha)\leq \prod_{j=1}^n t_j(\alpha_j n)^{1/n}
\quad(j=1,\ldots,n;\,\alpha\in\Nn).
\end{equation}
Consider the Fourier transforms
\begin{equation*}
\widehat\mu_k(\lambda)=\int_\Rn
e^{i(\lambda,x)}\,d\mu_k(x)\quad(k=1,2;\,\lambda\in\Rn).
\end{equation*}
Then $\widehat\mu_1$ and $\widehat\mu_2$ are of class $C^\infty$ on $\Rn$ with
derivatives
\begin{equation}\label{eq:derivatives1}
\frac{\partial^\alpha \widehat\mu_k}{\partial\lambda^\alpha}(\lambda)=\int_\Rn
i^{\vert\alpha\vert}x^\alpha
e^{i(\lambda,x)}\,d\mu_k(x)\quad(k=1,2;\,\alpha\in\Nn;\,\lambda\in\Rn).
\end{equation}
By assumption we therefore have
\begin{equation}\label{eq:equalderivatives}
\frac{\partial^\alpha \widehat\mu_1}{\partial \lambda^\alpha}(0)=\frac{\partial^\alpha \widehat\mu_2}{\partial
\lambda^\alpha}(0)\quad(\alpha\in\Nn).
\end{equation}
From \eqref{eq:derivatives1} we see that
\begin{equation*}
\frac{1}{2}\left\vert \frac{\partial^\alpha(\widehat\mu_1-\widehat\mu_2)}{\partial\lambda^\alpha}(\lambda)
\right\vert\leq t(\alpha)\quad(\alpha\in\Nn;\,\lambda\in\Rn),
\end{equation*}
and then combination with \eqref{eq:Holder} yields
\begin{equation*}
\frac{1}{2}\left\vert
\frac{\partial^\alpha(\widehat\mu_1-\widehat\mu_2)}{\partial\lambda^\alpha}
(\lambda)\right\vert\leq\prod_{j=1}^n t_j(\alpha_j
n)^{1/n}\quad(\alpha\in\Nn;\,\lambda\in\Rn).
\end{equation*}
We claim that the non-negative sequences $\{t_j(mn)^{1/n}\}_{m=0}^\infty$ ($j=1,\ldots,n$)
satisfy the hypotheses of Theorem~\ref{thm:DCtheorem}, i.e.\ we claim that
\begin{equation}\label{eq:divergenceCarleman1}
\sum_{m=1}^\infty \frac{1}{t_j(mn)^{1/mn}}=\infty.
\end{equation} To see this we fix $j$. If $\vert x_j\vert=0$ almost
everywhere ($\nu$) then \eqref{eq:divergenceCarleman1} is obvious. If $\vert
x_j\vert$ is not $\nu$-almost everywhere equal to zero, we define the (then
finite valued) non-negative sequence $\{h_j(m)\}_{m=1}^\infty$ by
\begin{equation*}
h_j(m)=t_j(m)^{-1/m}\quad(m=1,2,\ldots).
\end{equation*}
By \eqref{eq:increasing} the sequence $\{h_j(m)\}_{m=1}^\infty$ is
non-increasing. For $m$ even we have $s_j(m)=t_j(m)$; the hypothesis
\eqref{eq:Carleman} therefore translates as $\sum_{m=1}^\infty h_j(2m)=\infty$.
Lemma~\ref{lem:elementary} then implies that $\sum_{m=1}^\infty h_j(mn)=\infty$,
which is \eqref{eq:divergenceCarleman1}. This establishes the claim.

To conclude the proof of the determinacy we note that by
\eqref{eq:equalderivatives} all derivatives of
$\frac{1}{2}(\widehat\mu_1-\widehat\mu_2)$ vanish at $0$. Therefore Theorem
\ref{thm:DCtheorem} now shows that $\widehat \mu_1=\widehat\mu_2$, implying that
$\mu_1=\mu_2$, as was to be proved.

We turn to the density statements in $L_p(\Rn,\mu)$ for $1\leq p<\infty$. Fix
such $p$ and let $1<q\leq\infty$ be the conjugate exponent.

We treat the polynomials first. Suppose $f\in L_q(\Rn,\mu)$ is such that
\begin{equation}
\int_\Rn P(x) f(x)\,d\mu(x)=0
\end{equation}
for all polynomials $P$. We need to prove that $f=0$ a.e.\ ($\mu$). Define the
complex Borel measure $\xi_f$ on $\Rn$ by
\begin{equation}\label{eq:newmeasure}
\xi_f(E)=\int_E f(x)\,d\mu(x)
\end{equation}
for Borel sets $E$. Consider the Fourier transform
\begin{equation}\label{eq:Fourierxif}
\widehat\xi_f(\lambda)=\int_\Rn e^{i(\lambda,x)}\,d\xi_f(x)=\int_\Rn e^{i(\lambda,x)}
f(x)\,d\mu(x).
\end{equation}
Then $\widehat\xi_f$ is of class $C^\infty$ on $\Rn$ with derivatives
\begin{equation}\label{eq:derivatives2}
\frac{\partial^\alpha \widehat\xi_f}{\partial\lambda^\alpha}(\lambda)=\int_\Rn
i^{\vert\alpha\vert}x^\alpha
e^{i(\lambda,x)}f(x)\,d\mu(x)\quad(\alpha\in\Nn;\,\lambda\in\Rn).
\end{equation}
By assumption we therefore have
\begin{equation}\label{eq:derivativesvanish}
\frac{\partial^\alpha \widehat\xi_f}{\partial\lambda^\alpha}(0)=0\quad(\alpha\in\Nn).
\end{equation}
For $\alpha\in\Nn$ and $\lambda\in\Rn$ we have the following estimate, as a
consequence of \eqref{eq:derivatives2} and the generalized H\"older inequality
(the norms refer to $\mu$):
\begin{align*}
\left\vert \frac{\partial^\alpha
\widehat\xi_f}{\partial\lambda^\alpha}(\lambda)\right\vert &\leq\left\Vert\, \vert
x\vert^\alpha
\vert f\vert\,\right\Vert_1\\
&\leq \left\Vert f\right\Vert_q\left\Vert\,\vert x\vert^\alpha\right\Vert_p \\
&=\left\Vert f\right\Vert_q\Vert \prod_{j=1}^n
\vert x_j\vert^{\alpha_j p}\Vert_1^{1/p}\\ &\leq \left\Vert f\right\Vert_q
\prod_{j=1}^n
\left\Vert\,\vert x_j\vert^{\alpha_j p}\right\Vert_n^{1/p} \\&=\left\Vert f\right\Vert_q\prod_{j=1}^n t_j(\alpha_j p n
)^{1/np}.
\end{align*}
Here we have used \eqref{eq:marginalsequences} for the definition of $t_j(s)$
($j=1,\ldots,n;\,s\geq 0$), which is correct since we already know that
$\nu=\mu$. We claim that the non-negative sequences
$\{t_j(mpn)^{1/np}\}_{m=0}^\infty$ ($j=1,\ldots,n$) satisfy the hypotheses of Theorem
\ref{thm:DCtheorem}, i.e.\ we claim that
\begin{equation}\label{eq:divergenceCarleman2}
\sum_{m=1}^\infty \frac{1}{t_j(mpn)^{1/mnp}}=\infty.
\end{equation} To see this we again fix $j$. If $\vert x_j\vert=0$ almost
everywhere ($\mu$) then \eqref{eq:divergenceCarleman2} is again obvious. If
$\vert x_j\vert$ is not $\mu$-almost everywhere equal to zero, then we note that
\eqref{eq:increasing} implies that
\begin{equation}\label{eq:divergenceCarleman3}
\sum_{m=1}^\infty \frac{1}{t_j(mpn)^{1/mnp}}\geq \sum_{m=1}^\infty
\frac{1}{t_j(m([p]+1)n)^{1/m([p]+1)n}},
\end{equation}
where $[p]$ is the largest integer not exceeding $p$. In the notation as in the
proof of the determinacy, the right hand side of \eqref{eq:divergenceCarleman3}
is $\sum_{m=1}^\infty h_j(m([p]+1)n)$. Again Lemma~\ref{lem:elementary} then
implies that this series is divergent since $\sum_{m=1}^\infty h_j(2m)$
diverges, thus establishing the claim.

In view of \eqref{eq:derivativesvanish} we now conclude from Theorem
\ref{thm:DCtheorem} that $\widehat\xi_f$=0, implying $\xi_f=0$ and finally that
$f=0$ a.e.\ ($\mu$), as was to be proved.

Finally, let us prove the density of $\trigspan$ for a subset $S$ of $\Rn$ such that $\overline
S$ has non-empty interior. Assume that $a\in\Rn$ is an interior point of $\overline S$ and
suppose $f\in L_q(\Rn,\mu)$ vanishes on $\trigspan$. Consider the Fourier transform
\begin{equation}\label{eq:trigonometrics}
\widehat{e_{ia}f}(\lambda)=\int_\Rn e^{i(\lambda+a,x)}f(x)\,d\mu(x)\quad (\lambda\in
\Rn).
\end{equation}
Then $\widehat{e_{ia}f}$ is of class $C^\infty$ and the assumption on $f$ implies that
$\widehat{e_{ia}f}$ is identically zero on a neighborhood of $0\in\Rn$. Evidently all
derivatives of $\widehat{e_{ia}f}$ then vanish at $0$, which shows that $e_{ia}f$ vanishes on
the polynomials. Since we had already shown that these are dense in $L_p(\Rn,\mu)$ we conclude
that $f=0$ a.e.\ ($\mu$), as was to be proved.
\end{proof}

We comment on the relation between Theorem~\ref{thm:Carleman} and the
literature.

The fact that the divergence of the series in Theorem~\ref{thm:Carleman} is sufficient for the
determinacy of the measure can already be found in \cite{Nussbaum}, where a combination of
quasi-analytic methods and a Hilbert space approach is used.

Theorem~\ref{thm:Carleman} is also related to the following result \cite[p.~21]{ShohatTamarkin}:
if $\mu\in\measures$ and if we have $\sum_{m=1}^\infty 1/\sqrt[2m]{\lambda(2m)}=\infty$, where
\begin{equation}
\lambda(m)=\int_\Rn \sum_{j=1}^n x_j^m\,d\mu(x)\quad(m=0,1,2,\ldots),
\end{equation}
then $\mu$ is determinate. It is under this condition on the $\lambda(2m)$ even true
\cite{Berg1,Berg2} that the polynomials are dense in $L_p(\Rn,\mu)$ for all $1\leq p<\infty$, a
property which is stronger than determinacy of the measure. As to this last implication, it was
first proved in \cite{BergChristensen1} in the one dimensional case that for $\mu\in\measures$
the density of the polynomials in $L_p(\R,\mu)$ for some finite $p>2$ implies that $\mu$ is
determinate, and this result was later generalized to arbitrary dimension in \cite{Fuglede}.

The determinacy and polynomial density under the condition on the $\lambda(2m)$ also follow from
Theorem~\ref{thm:Carleman}. Indeed, taking the standard basis in Theorem~\ref{thm:Carleman} one
obviously has $s_j(2m)\leq
\lambda(2m)$. The divergence of the series for $\lambda(2m)$ therefore implies the divergence of
the series for all $s_j(m)$, so that the conclusions of Theorem~\ref{thm:Carleman} on
determinacy and density hold.

Conversely, the special case $\R^n=\R$ (see \cite{BergChristensen2}) of the results on
determinacy and polynomial density as quoted above can be taken as a starting point to derive
Theorem~\ref{thm:Carleman}, albeit in a more indirect fashion than in the present paper. Indeed,
assuming that the basis in Theorem~\ref{thm:Carleman} is the standard basis, one concludes from
this one dimensional starting point that all marginal distributions of the measure in Theorem
\ref{thm:Carleman} are determinate and that the polynomials are dense in all $L_p$-spaces
($1\leq p<\infty$) associated with these marginal distributions. The results in \cite{Petersen}
then imply that the analogous two statements hold for the measure itself. The additional density
of the trigonometric functions then also follows. Indeed, the conclusion of the proof of Theorem
\ref{thm:Carleman} shows in fact that for $\mu\in\measures$ and $1\leq p<\infty$, the density of
the polynomials in $L_p(\Rn,\mu)$ implies the density in $L_p(\Rn,\mu)$ of $\trigspan$ for all
somewhere dense subsets $S$ of $\Rn$. The author is indebted to Christian Berg for communicating
this last result and its proof.

\section{Quasi-analytic weights}\label{sec:qaw}

This section is a preparation for Sections~\ref{sec:Hamburger} and \ref{sec:Stieltjes}. We
define quasi-analytic weights, mention the relevant properties and the classification and give practical examples.

The notion of quasi-analytic weight is a delicate one, which is studied systematically in
\cite{deJeu}. Simple and intuitive properties are sometimes not immediately obvious and may
require an argument. We therefore refer for proofs to \cite{deJeu}, from which all results in
this section are taken.

A weight on $\Rn$ is an arbitrary bounded non-negative function on $\Rn$. We emphasize that we
assume no regularity.

\begin{definition}\label{def:qaw}
Let $w$ be a weight on $\Rn$. Suppose $\{v_1,\ldots,v_n\}$ is a basis of $\Rn$.
If
\begin{equation*}
\sum_{m=1}^\infty \frac{1}{{\Vert (v_j,x)^m w(x)\Vert_\infty^{1/m}}}=\infty
\end{equation*}
for $j=1,\ldots,n$ then $w$ is \emph{quasi-analytic with respect to $\{v_1,\ldots,v_n\}$}. A weight is \emph{standard quasi-analytic} if it is quasi-analytic with respect to the standard basis of $\Rn$. A weight is
\emph{quasi-analytic} if it is quasi-analytic with respect to some basis.
\end{definition}

The terminology ``quasi-analytic'' refers not to regularity of the weight itself, which might
e.g.\ even fail to be Lebesgue measurable. The reason then for this terminology lies---as has
become customary---in the fact that certain crucial associated functions have the quasi-analytic
property, by which is meant that one can conclude that such an associated function is actually
identically zero once one has established that the function and all its derivatives vanish at
one fixed point. These associated functions thus share this property with analytic functions in
function theory, which explains the nomenclature. In our case, the reader may verify that it is
indeed the divergence of the series in Definition~\ref{def:qaw} that validates the application
of Theorem~\ref{thm:DCtheorem} on the quasi-analytic property in the proof of the basic
Theorem~\ref{thm:Carleman}.

A weight which vanishes outside a compact set is quasi-analytic with respect to
all bases. By a small computation, the same holds for a weight of type
$\exp(-\epsilon \Vert x\Vert)$ with $\epsilon>0$. The set of quasi-analytic
weights is invariant under the group of affine automorphisms of $\Rn$ and under
multiplication with non-negative constants.

Let $w$ be a weight, quasi-analytic with respect to the basis
$\{v_1,\ldots,v_n\}$. Then $\lim_{\Vert x\Vert\rightarrow\infty}\Vert
x\Vert^dw(x)=0$ for all $d\geq 0$, i.e.\ $w$ is rapidly decreasing. A minorant
of $w$ outside a compact set is again quasi-analytic with respect to $\{v_1,\ldots,v_n\}$.
One can prove that $w$ always has a pointwise majorant of class $C^\infty$ which
is strictly positive and quasi-analytic with respect to $\{v_1,\ldots,v_n\}$.
In one dimension, such a majorant can in addition be required to be even and strictly decreasing on $[0,\infty)$.

There are several closely related ways of characterizing quasi-analytic weights
other than by Definition~\ref{def:qaw}, which is a technically convenient
characterization but not a very practical one to verify. The formulation in the following
paragraphs seems to fit most applications. For additional material the reader is
referred to \cite{deJeu}.

A weight $w$ on $\Rn$ is quasi-analytic if and only if there exists an affine
automorphism $A$ of $\Rn$ and quasi-analytic weights $w_j$ ($j=1,\ldots,n$) on
$\R$ such that
\begin{equation}\label{eq:classndim}
w(Ax)\leq \prod_{j=1}^m w_j(x_j)
\end{equation}
for all $x\in\Rn$. More precisely, if $w$ is quasi-analytic with respect to the basis
$\{v_1,\ldots,v_n\}$ of $\Rn$, then quasi-analytic weights $w_j$ on $\R$ satisfying
\eqref{eq:classndim} exist for any $A$ with linear component $A_0$ defined by
$A_0^t v_j=e_j\,(j=1,\ldots,n)$; here $A_0^t$ is the transpose of $A_0$ with
respect to the standard inner product on $\Rn$. Conversely, if
\eqref{eq:classndim} holds for some $A$ and quasi-analytic weights $w_j$ on $\R$, and if $A_0$ is
the linear component of $A$, then $w$ is quasi-analytic with respect to the
basis $\{v_1,\ldots,v_n\}$ of $\Rn$ defined by $A_0^t v_j=e_j\,(j=1,\ldots,n)$.

The matter has now been reduced to $\R$. As a first equivalent characterization
on the real line, a weight $w$ on $\R$ is quasi-analytic if and only if there
exist $R>0$, $C\geq 0$ and a non-decreasing function
$\rho:(R,\infty)\mapsto\R_{\geq 0}$ of class $C^1$ such that
\begin{equation}\label{eq:divergentintegral}
\int_{R}^\infty\frac{\rho(s)}{s^2}\,ds=\infty
\end{equation}
and
\begin{equation*}
w(x)\leq C \exp\left(-\int_{R}^{\vert x\vert}\frac{\rho(s)}{s}\right)\,ds
\end{equation*}
if $\vert x\vert\geq R$.

As a second and closely related equivalent characterization on the real line, a
weight on $\R$ is quasi-analytic if and only if there exists a weight $\wtilde$
on $\R$ and $R>0$ such that $w(t)\leq
\wtilde(t)$ and $\wtilde(t)=\wtilde(-t)>0$ both hold for $\vert t\vert>R$,
such that $s\mapsto
-\log\wtilde (e^s)$ is convex on $(\log R,\infty)$ and such that
\begin{equation*}
\int_R^\infty\frac{\log\wtilde(t)}{1+t^2}=-\infty.
\end{equation*}
Weights such as $\wtilde$ are classical and figure e.g.\ in the Bernstein problem \cite{Koosis}.
The connection between these classical weights and the one dimensional version of Definition
\ref{def:qaw} seems to have gone largely unnoticed, although some ingredients can be found in
\cite[proof of Theorem 2]{Lin} under additional regularity conditions on the weight.

If $w$ is a quasi-analytic weight on $\R$, then one obtains a quasi-analytic
weight $w^\prime$ on $\Rn$ by putting $w^\prime(x)=w(\Vert x\Vert)$ for
$x\in\Rn$. Such $w^\prime$ is then quasi-analytic with respect to all bases of
$\Rn$. All minorants of $w^\prime$ outside a compact set are then again quasi-analytic
with respect to all bases of $\Rn$. The first alternative characterization of
quasi-analytic weights on $\R$, combined with this radial extension procedure
thus yields the following.

\begin{proposition}\label{prop:qawradial}
Suppose $R>0$ and a non-decreasing function $\rho:(R,\infty)\mapsto\R_{\geq 0}$
of class $C^1$ are such that
\begin{equation*}
\int_{R}^\infty\frac{\rho(s)}{s^2}\,ds=\infty.
\end{equation*}
If $w$ is a weight such that
\begin{equation*}
w(x)\leq C \exp\left(-\int_{R}^{\Vert x\Vert}\frac{\rho(s)}{s}\,ds\right)
\end{equation*}
whenever $\Vert x\Vert\geq R$, then $w$ is a weight on $\Rn$ which is
quasi-analytic with respect to all bases of $\Rn$.
\end{proposition}

The following result in terms of elementary functions is based on the second alternative
characterization of quasi-analytic weights on $\R$, combined with the radial extension
procedure.

\begin{proposition}\label{prop:examplesqaw}
Define repeated logarithms by $\log_0 t=t$ and, inductively, for $j\geq 1$, by
$\log_j t=\log (\log_{j-1}t)$, where $t$ is assumed to be sufficiently large for
the definition to be meaningful in the real context. For $j=0,1,2,\ldots$ let
$a_j>0$ and let $p_j\in\R$ be such that $p_j=0$ for all sufficiently large $j$.
Put $j_0=\min\{j=0,1,2,\ldots\mid p_j\neq 1\}$. Let $C>0$ and suppose
$w:\Rn\mapsto\R_{\geq 0}$ is bounded.

Then, if $p_{j_0}<1$ and if
\begin{equation*}
w(x)\leq C\exp\left(-\Vert x\Vert^2
\left(\prod_{j=0}^\infty\log_{j}^{p_j}a_j\Vert x\Vert\right)^{-1}\right)
\end{equation*}
for all sufficiently large $\Vert x\Vert$, $w$ is a weight on $\Rn$ which is
quasi-analytic with respect to all bases of $\Rn$.
\end{proposition}

Note the occurrence of $\log_0$ (i.e.\ of the identity) in the Proposition, which permits a
uniform formulation.

Thus, to give explicit examples, a weight on $\Rn$ is quasi-analytic if it is for all
sufficiently large $\Vert x\Vert$ majorized by one of the expressions
\begin{align*}
&C\exp\left(-\frac{\Vert x\Vert^{1-\nu}}{a_0}\right)\quad ,\\
&C\exp\left(-\frac{\Vert x\Vert}{a_0 \left(\log a_1\Vert
x\Vert\right)^{1+\nu}}\right)\quad ,\\&C\exp\left(-\frac{\Vert x\Vert}{a_0 \log
a_1\Vert x\Vert
\left(\log\log a_2\Vert x\Vert\right)^{1+\nu}}\right)\quad ,\\
&\ldots
\end{align*}
for some $C, a_0,a_1,a_2,\ldots>0$ and $\nu\leq 0$. The case $\nu=0$ yields a
sequence of families of quasi-analytic weights, each consisting of weights that
are negligible at infinity compared with any member of the succeeding family.

Explicit non-radial examples of standard quasi-analytic weights on $\Rn$ in
terms of elementary functions can be obtained as tensor products of
quasi-analytic weights on $\R$ taken from Proposition~\ref{prop:examplesqaw}.
All minorants of such tensor products outside a compact set are then again standard
quasi-analytic weights on $\Rn$. Insertion of an affine automorphism of $\Rn$ in
the argument of the weight yields additional quasi-analytic weights.

For the sake of completeness we mention that we have the following negative
criterion for a weight to be quasi-analytic: if $w$ is a quasi-analytic weight
on $\Rn$ and if $x,y\in\Rn$ are such that $y\neq 0$ and such that $t\mapsto
w(x+ty)$ is Lebesgue measurable on $\R$, then for all $R>0$ we have
\begin{equation*}
\int_R^\infty\frac{\log w(x+ty)}{1+t^2}\,dt=\int_{-\infty}^{-R}\frac{\log w(x+ty)}{1+t^2}\,dt=-\infty.
\end{equation*}
This shows that Proposition~\ref{prop:examplesqaw} is sharp in the sense that
the corresponding statement for $p_{j_0}>1$ does not hold.

The set of quasi-analytic weights has some interesting characteristics. Contrary to what the
explicit examples above suggest, this set is not closed under addition. More precisely, one can
construct weights $w_1$ and $w_2$ on $\Rn$, each of which is quasi-analytic with respect to all
bases, but such that $w_1+w_2$ is not quasi-analytic with respect to any basis. One can also
construct weights which are quasi-analytic with respect to just \emph{one} basis (up to
scaling); in Section~\ref{sec:closing} we will make some tentative remarks on a possible
parallel of this phenomenon for measures. For $n\geq 2$ it implies that such quasi-analytic
weights on $\Rn$ are not minorants outside a compact set of quasi-analytic weights as obtained
from the radial extension procedure.

\section{Integral criteria for determinacy}\label{sec:Hamburger}

We will now combine the results in Sections~\ref{sec:Carleman} and~\ref{sec:qaw}.

\begin{theorem}[First main theorem]\label{thm:maintheoremHamburger}
Let $\mu$ be a positive Borel measure on $\Rn$ such that
\begin{equation*}
\int_\Rn w(x)^{-1}\,d\mu<\infty
\end{equation*}
for some measurable quasi-analytic weight. Then $\mu$ is determinate.
Furthermore, the polynomials and $\trigspan$ are then dense in $L_p(\Rn,\mu)$
for all $1\leq p<\infty$ and for every subset $S$ of $\Rn$ which is somewhere
dense.
\end{theorem}

Note that since quasi-analytic weights are rapidly decreasing, the measure in the theorem is automatically in $\measures$.

\begin{proof}{Proof}
Suppose $w$ is quasi-analytic with respect to the basis $\{v_1,\ldots,v_n\}$. We may assume that
$w$ is strictly positive: if necessary we can replace $w$ by a strictly positive measurable (say
smooth) majorant which is quasi-analytic with respect to $\{v_1,\ldots,v_n\}$. We may also
assume that $\Vert w\Vert_\infty=1$.

In the notation of Theorem~\ref{thm:Carleman} we then have for $j=1,\ldots,n$
and $m=0,1,2,\ldots$:
\begin{align}\label{eq:estimates}
\begin{split}
s_j(2m)&=\int_\Rn (v_j,x)^{2m} \,d\mu(x)\\ &=\int_\Rn (v_j,x)^{2m} w(x) w(x)^{-1}\,d\mu(x)\\
&\leq \Vert (v_j,x)^{2m} w(x)\Vert_\infty\int_\Rn w(x)^{-1}\,d\mu(x).
\end{split}
\end{align}
Now the sequences $\{\Vert (v_j,x)^m w(x)\Vert_\infty^{1/m}\}_{m=1}^\infty$ ($j=1,\ldots,n$) are
easily seen to be non-decreasing, as a consequence of the normalization $\Vert w\Vert_\infty=1$.
The quasi-analyticity of $w$ with respect to $\{v_1,\ldots,v_n\}$ therefore implies by Lemma
\ref{lem:elementary} that
\begin{equation*}
\sum_{m=1}^\infty \frac{1}{{\Vert (v_j,x)^{2m}
w(x)\Vert_\infty^{1/2m}}}=\infty\quad(j=1,\ldots,n).
\end{equation*}
This divergence implies, in view of the estimate in \eqref{eq:estimates}, that the hypotheses of
Theorem~\ref{thm:Carleman} are satisfied, which concludes the proof.
\end{proof}

Theorem~\ref{thm:maintheoremHamburger} can also be found in \cite{deJeu}, where the
density part is seen to be a consequence of more general considerations on the
closure of modules over the polynomials and trigonometric functions in
topological vector spaces. In [loc.cit.] the determinacy of the measure is then
concluded from \cite{Fuglede} since there exists $p>2$ such that the polynomials are dense in the associated $L_p$-space. This way of deriving Theorem
\ref{thm:maintheoremHamburger} is considerably more involved than the present proof.

The combination of Theorem~\ref{thm:maintheoremHamburger} with the results on quasi-analytic
weights in Section~\ref{sec:qaw} now yields various integral criteria for determinacy and
density, as mentioned in the introduction. Variation in these criteria, in particular variation
in the degree of regularity of the weights involved, is possible in view of the various ways in
which quasi-analytic weights can be characterized. In the criteria as described in this section,
all integrands are of class $C^2$ outside a compact set. For cases where this is too stringent
the reader is referred to \cite{deJeu}.

A non-radial criterion is the following.

\begin{theorem}\label{thm:criteriumHamburger}
For $j=1,\ldots,n$ let $R_j>0$ and a non-decreasing function
$\rho_j:(R_j,\infty)\mapsto\R_{\geq 0}$ of class $C^1$ be such that
\begin{equation*}
\int_{R_j}^\infty\frac{\rho_j(s)}{s^2}\,ds=\infty.
\end{equation*}
Define $f_j:\R\mapsto\R_{\geq 0}$ by
\begin{equation*}
f_j(x)= \exp{\left(\int_{R_j}^{\vert x\vert}\frac{\rho_j(s)}{s}\,ds\right)}
\end{equation*}
for $\vert x\vert> R_j$ and by $f_j(x)=1$ for $\vert x\vert\leq R_j$. Let $A$ be
an affine automorphism of $\Rn$. If $\mu$ is a positive Borel measure on $\Rn$
such that
\begin{equation*}
\int_\Rn \prod_{j=1}^n f_j((Ax)_j)\,d\mu(x)<\infty,
\end{equation*}
then $\mu$ is determinate. Furthermore, the polynomials and $\trigspan$ are then
dense in $L_p(\Rn,\mu)$ for all $1\leq p<\infty$ and for every subset $S$ of
$\Rn$ which is somewhere dense.
\end{theorem}

\begin{proof}{Proof}
From the first alternative characterization of quasi-analytic weights on $\R$, as given in
Section~\ref{sec:qaw} we see that the weights $1/f_j$ are all quasi-analytic weights on $\R$.
Their tensor product is then a quasi-analytic weight on $\Rn$ and then the same holds for the
image of this tensor product under an element of the affine group. We now apply Theorem
\ref{thm:maintheoremHamburger}.
\end{proof}

We now specialize to the case of radial integrands. Theorem~\ref{thm:introthm} evidently follows
from Theorem~\ref{thm:maintheoremHamburger} and Proposition~\ref{prop:qawradial}. In addition, the
combination of Proposition~\ref{prop:examplesqaw} and Theorem~\ref{thm:maintheoremHamburger} implies the
following. As with Proposition~\ref{prop:examplesqaw}, note the occurrence of $\log_0$, i.e.\ of
the identity.

\begin{theorem}\label{thm:detexplicit}
Define repeated logarithms $\log_j$ ($j=0,1,2,\ldots$) as in Proposition
\ref{prop:examplesqaw}. For $j=0,1,2,\ldots$ let $a_j>0$ and let $p_j\in\R$ be
such that $p_j=0$ for all sufficiently large $j$. Put
$j_0=\min\{j=0,1,2,\ldots\mid p_j\neq 1\}$ and assume $p_{j_0}<1$.

Let $\mu$ be a positive Borel measure on $\Rn$. If

\begin{equation*}
\int_{\Vert x\Vert\geq R}\exp\left(\Vert x\Vert^2\left(\prod_{j=0}^\infty\log_{j}^{p_j}a_j\Vert
x\Vert\right)^{-1}\right)\,d\mu<\infty,
\end{equation*}
for some $R\geq 0$ which is sufficiently large to ensure that the integrand is
defined, then $\mu$ is determinate. Furthermore, the polynomials and $\trigspan$
are then dense in $L_p(\Rn,\mu)$ for all $1\leq p<\infty$ and for every subset
$S$ of $\Rn$ which is somewhere dense.
\end{theorem}

As explicit examples, if one of the functions (tacitly assumed to be equal to $1$ on a
sufficiently large compact set)
\begin{align*}
&\exp\left(\frac{\Vert x\Vert^{1-\nu}}{a_0}\right)\quad ,\\ &\exp\left(\frac{\Vert x\Vert}{a_0
\left(\log a_1\Vert x\Vert\right)^{1+\nu}}\right)\quad ,\\&\exp\left(\frac{\Vert
x\Vert}{a_0 \log a_1\Vert x\Vert
\left(\log\log a_2\Vert x\Vert\right)^{1+\nu}}\right)\quad ,\\
&\ldots
\end{align*}
has finite integral under $\mu$ for some $a_0,a_1,a_2,\ldots>0$ and $\nu\leq 0$, then the
conclusions in Theorem~\ref{thm:detexplicit} hold. The classical condition of the integrability
of $\exp(\epsilon \Vert x\Vert)$ for some $\epsilon>0$ can be weakened quite substantially.

To conclude we mention that explicit non-radial reciprocals of quasi-analytic weights can be
obtained in terms of elementary functions by taking the tensor product of the reciprocals of the
one dimensional versions of the majorants in Proposition~\ref{prop:examplesqaw}, when these
majorants are in addition defined to be equal to $1$ on a sufficiently large compact subset of
$\R$.

\section{Determinacy in the sense of Stieltjes}\label{sec:Stieltjes}

In this section we are concerned with determinacy in the sense of Stieltjes, i.e.\ with
$C$-determinacy as in Definition~\ref{def:cdet}. Analogously to the one dimensional case, The
Carleman criterion in Theorem~\ref{thm:Carleman} implies a similar sufficient condition for
$C$-determinacy. When combined with the results on quasi-analytic weights again, we obtain
integral criteria for $C$-determinacy. At the end of the section we discuss a condition (which
is satisfied for absolutely continuous measures) enabling one to conclude that the measure is
not just $C$-determinate, but in fact determinate.

\begin{theorem}[Carleman criterion for $C$-determinacy]\label{thm:CarlemanStieltjes}
Let $\{v_1,\ldots,v_n\}$ be a basis of $\Rn$ and let $C=\R_{\geq 0}\cdot v_1+\ldots\R_{\geq
0}\cdot v_n$ be the corresponding positive convex cone. Define the dual basis
$\{v_1^\prime,\ldots,v_n^\prime\}$ by $(v_i^\prime,v_j)=\delta_{ij}$ $(i,j=1,\ldots,n)$.

Let $\mu\in\measures$ be supported in $C$. For $j=1,\ldots,n$ and $m=0,1,2,\ldots$ define
\begin{equation*}
s_j(m)=\int_C (v_j^\prime,x)^m\,d\mu(x),
\end{equation*}
and suppose that each of the sequences $\{s_j(m)\}_{m=1}^\infty$ ($j=1,\ldots,n$) satisfies
\begin{equation}
\sum_{m=1}^\infty\frac{1}{s_j(m)^{1/2m}}=\infty.
\end{equation}

Then $\mu$ is $C$-determinate.
\end{theorem}

Note that the $s_j(m)$ are defined in terms of distinguished coordinates on $C$, namely those
corresponding to extremal generators of $C$.

\begin{proof}{Proof}
The proof generalizes the well known proof in one dimension. Let $\measuresC$ be the measures in
$\measures$ which are supported in $C$. As a first preparation, define $\phi:\Rn\mapsto\Rn$ by
\begin{equation*}
\phi(x)=
\begin{cases}
x\quad &\text{if }x\notin C;\\
\sum_{j=1}^n \sqrt{x_j}v_j\quad&\text{if }x=\sum_{j=1}^n x_j v_j,\,x_j\geq 0\,\,(j=1,\ldots,n).
\end{cases}
\end{equation*}
For $\xi\in\measures$, define $\xi_\phi\in\measures$ by putting $\xi_\phi(A)=\xi(\phi^{-1}(A))$
for a Borel set $A$. The assignment $\xi\mapsto\xi_\phi$ defines an injective map from
$\measures$ to $\measures$ which maps $\measuresC$ into itself, and
\begin{equation}\label{eq:trafoformula}
\int_\Rn P(x)\,d\xi_\phi(x)=\int_\Rn (P\circ\phi)(x)\,d\xi(x)
\end{equation}
for all polynomials $P$ and all $\xi\in\measures$.

As a second preparation, let $G$ be the group of linear isomorphisms of $\Rn$ having $2^n$
elements, corresponding to all possible sign changes in the coordinates with respect to the
basis $\{v_1,\ldots,v_n\}$. For $\xi\in\measures$ and $g\in G$ define $g\cdot
\xi\in\measures$ by putting $(g\cdot\xi)(A)=\xi(g^{-1}(A))$ for a Borel set $A$ and let
$\overline\xi=2^{-n}\sum_{g\in G} g\cdot
\xi$. The averaging map $\xi\mapsto\overline\xi$ is not injective as a map from $\measures$ to
itself, but it \emph{is} injective as a map from $\measuresC$ to $\measures$. To see this, let
$J\subset\{1,\ldots,n\}$ have cardinality $\vert J\vert$ and define $C_J=\{x\in C\mid x_j=0
\Leftrightarrow j\in J\}$. Then for $\xi\in\measuresC$ we have $\overline\xi(A_J)= 2^{\vert
J\vert-n}\xi(A_J)$ for any Borel subset $A_J$ of $C_J$. Thus the restriction of $\xi$ to $C_J$
can be retrieved from $\overline\xi$. Since the $C_J$ form a disjoint covering of $C$, we see
that $\xi\in\measuresC$ can be reconstructed from $\overline\xi$, as claimed.

Furthermore, if $\xi\in\measures$, then
\begin{equation}\label{eq:invariantintegral}
\int_\Rn (P\circ g)(x)\,d\,\overline\xi(x)=\int_\Rn P(x)\,d\,\overline\xi(x)
\end{equation}
for all polynomials $P$, so that
\begin{equation}\label{eq:integralzero}
\int_\Rn \prod_{j=1}^n(v_j^\prime,x)^{e_j} \,d\overline\xi(x)=0
\end{equation}
if the $e_j$ are non-negative integers, at least one of which is odd. This follows from \eqref{eq:invariantintegral} by choosing an element in $G$ which sends the integrand in the left hand side of \eqref{eq:integralzero} to its negative. On the other hand, if
$\xi\in\measures$ and if the $e_j$ are all even non-negative integers, then
\begin{equation}\label{eq:integralequal}
\int_\Rn \prod_{j=1}^n(v_j^\prime,x)^{e_j} \,d\overline\xi(x)=\int_\Rn \prod_{j=1}^n(v_j^\prime,x)^{e_j}
\,d\xi(x)
\end{equation}
as a consequence of the invariance of the integrand under $G$.

Combining \eqref{eq:trafoformula}, \eqref{eq:integralzero} and \eqref{eq:integralequal}, we
conclude that for $\xi\in\measuresC$
\begin{equation}\label{eq:momentsrelation}
\int_\Rn \prod_{j=1}^n(v_j^\prime,x)^{e_j} \,d\overline{\xi_\phi}(x)=\int_C \prod_{j=1}^n(v_j^\prime,x)^{e_j/2}
\,d\xi(x)
\end{equation}
if the $e_j$ are all even non-negative integers, whereas the integral is zero if the $e_j$ are
non-negative integers, at least one of which is odd.

Turning to the theorem, we first of all note that $\overline{\mu_\phi}$ satisfies the conditions of Theorem~\ref{thm:Carleman} as a consequence of \eqref{eq:momentsrelation}, so $\overline{\mu_\phi}$ is determinate as a measure on $\Rn$. Suppose then that $\nu\in\measuresC$ yields the same integrals for all polynomials as $\mu$. Then \eqref{eq:momentsrelation} implies that $\overline{\mu_\phi}$ and $\overline{\nu_\phi}$ also have the same integrals for all
polynomials and we conclude that $\overline{\mu_\phi}=\overline{\nu_\phi}$. By the injectivity of
the maps as observed above, it first follows that $\mu_\phi=\nu_\phi$ and subsequently that
$\mu=\nu$.
\end{proof}

We will now combine this with the results in Section~\ref{sec:qaw}.

\begin{theorem}[Second main theorem]\label{thm:maintheoremStieltjes}
Let $\{v_1,\ldots,v_n\}$ be a basis of $\Rn$ and let $C=\R_{\geq 0}\cdot v_1+\ldots\R_{\geq
0}\cdot v_n$ be the corresponding positive convex cone. Define the dual basis
$\{v_1^\prime,\ldots,v_n^\prime\}$ by $(v_i^\prime,v_j)=\delta_{ij}$ $(i,j=1,\ldots,n)$. Let $w$
be a measurable weight on $\Rn$, quasi-analytic with respect to
$\{v_1^\prime,\ldots,v_n^\prime\}$. For $x=\sum_{j=1}^n x_j v_j\in C$ define
$\phi(x)=\sum_{j=1}^n \sqrt{x_j}v_j$. Let $\mu$ be a positive Borel measure on $\Rn$ which is
supported in $C$ and suppose that
\begin{equation*}
\int_C (w\circ\phi)(x)^{-1}\,d\mu(x)<\infty.
\end{equation*}
Then $\mu$ is $C$-determinate.
\end{theorem}

\begin{proof}{Proof}
The argument parallels the proof of Theorem~\ref{thm:maintheoremHamburger}. We may assume that
$w$ is strictly positive: if necessary we can replace $w$ by a strictly positive measurable (say
smooth) majorant which is quasi-analytic with respect to $\{v_1^\prime,\ldots,v_n^\prime\}$. We
may also assume that $\Vert w\Vert_\infty=1$.

In the notation of Theorem~\ref{thm:CarlemanStieltjes} we then have for $j=1,
\ldots,n$ and $m=0,1,2,\ldots$:
\begin{equation*}
s_j(m)\leq \left[\sup_{x\in C}\,(v_j^\prime,x)^m(w\circ\phi)(x)\right]\cdot \int_C
(w\circ\phi)(x)^{-1}\,d\mu(x).
\end{equation*}
Now
\begin{align*}
\begin{split}
\sup_{x\in C}\,(v_j^\prime,x)^m(w\circ\phi)(x)&=\sup_{x_1,\ldots,x_n\geq 0}
(v_j^\prime,\sum_{j=1}^n x_j v_j)^m w(\sum_{j=1}^n \sqrt{x_j} v_j)\\ &=\sup_{t_1,\ldots,t_n\geq
0} (v_j^\prime,\sum_{j=1}^n t_j^2 v_j)^m w(\sum_{j=1}^n t_j v_j)
\\&=\sup_{x\in C}\,(v_j^\prime,x)^{2m} w(x)\\&\leq \Vert (v_j^\prime,x)^{2m} w\Vert_\infty.
\end{split}
\end{align*}
As in the proof of Theorem~\ref{thm:maintheoremHamburger}, we conclude from these estimates and
Lemma~\ref{lem:elementary} that the hypotheses of Theorem~\ref{thm:CarlemanStieltjes} are
satisfied.
\end{proof}

The proof of the following theorem is left to the reader. It follows from Theorem~\ref{thm:maintheoremStieltjes}, using
the results in Section~\ref{sec:qaw} on quasi-analytic
weights, in a way similar to the proof of Theorem~\ref{thm:criteriumHamburger}.

\begin{theorem}\label{thm:tensorStieltjes}
Let $\{v_1,\ldots,v_n\}$ be a basis of $\Rn$ and let $C=\R_{\geq 0}\cdot v_1+\ldots\R_{\geq
0}\cdot v_n$ be the corresponding positive convex cone. Define the dual basis
$\{v_1^\prime,\ldots,v_n^\prime\}$ by $(v_i^\prime,v_j)=\delta_{ij}$ $(i,j=1,\ldots,n)$. For
$x=\sum_{j=1}^n x_j v_j\in C$ define $\phi(x)=\sum_{j=1}^n \sqrt{x_j}v_j$.

For $j=1,\ldots,n$ let $R_j>0$ and a non-decreasing function
$\rho_j:(R_j,\infty)\mapsto\R_{\geq 0}$ of class $C^1$ be such that
\begin{equation*}
\int_{R_j}^\infty\frac{\rho_j(s)}{s^2}\,ds=\infty.
\end{equation*}
Define $f_j:\R\mapsto\R_{\geq 0}$ by
\begin{equation*}
f_j(x)= \exp{\left(\int_{R_j}^{\vert x\vert}\frac{\rho_j(s)}{s}\,ds\right)}
\end{equation*}
for $\vert x\vert> R_j$ and by $f_j(x)=1$ for $\vert x\vert\leq R_j$.

If $\mu$ is a positive Borel measure on $\Rn$ which is supported in $C$, and if
\begin{equation*}
\int_C \prod_{j=1}^n f_j((v_j^\prime,\phi(x)))\,d\mu(x)<\infty,
\end{equation*}
then $\mu$ is $C$-determinate.
\end{theorem}

We turn to radial integrands, for which we will use the following practical result as a starting
point.

\begin{theorem}\label{thm:radialStieltjes}
Let $\{v_1,\ldots,v_n\}$ be a basis of $\Rn$ and let $C=\R_{\geq 0}\cdot v_1+\ldots\R_{\geq
0}\cdot v_n$ be the corresponding positive convex cone. Let $w$ be a measurable quasi-analytic
weight on the real line. Let $\mu$ be a positive Borel measure on $\Rn$ which is supported in
$C$ and suppose that
\begin{equation*}
\int_C w(\sqrt{\Vert x\Vert})^{-1}\,d\mu(x)<\infty.
\end{equation*}

Then $\mu$ is $C$-determinate.
\end{theorem}

\begin{proof}{Proof}
As with Theorem~\ref{thm:maintheoremStieltjes}, the proof parallels that of Theorem
\ref{thm:maintheoremHamburger}. We may assume that $w$ is strictly positive, by
replacing $w$ with a quasi-analytic majorant with this property if necessary. We may also assume
that $\Vert w\Vert_\infty=1$. In the notation of Theorem~\ref{thm:CarlemanStieltjes} we then
have for $j=1,
\ldots,n$ and $m=0,1,2,\ldots$:
\begin{equation*}
s_j(m)\leq \left[\sup_{x\in C}\,(v_j^\prime,x)^m w(\sqrt{\Vert x\Vert})\right]\cdot \int_C
w(\sqrt{\Vert x\Vert})^{-1}\,d\mu(x).
\end{equation*}
Now
\begin{align*}
\sup_{x\in C}\, (v_j^\prime,x)^m w(\sqrt{\Vert x\Vert})&\leq \Vert v_j^\prime\Vert^m \sup_{x\in C}\,\Vert x\Vert^{m} w(\sqrt{\Vert x\Vert})\\
&\leq\Vert v_j^\prime\Vert^m \sup_{t\in\R}\vert t^{2m} w(t)\vert.
\end{align*}
As in the proof of Theorem~\ref{thm:maintheoremHamburger}, we conclude from these estimates and
Lemma~\ref{lem:elementary} that the hypotheses of Theorem~\ref{thm:CarlemanStieltjes} are
satisfied.
\end{proof}

Theorem~\ref{thm:introthmStieltjes} is now obvious, given Theorem~\ref{thm:radialStieltjes} and Proposition~\ref{prop:qawradial}.
The combination of Proposition~\ref{prop:examplesqaw} and Theorem~\ref{thm:radialStieltjes} implies the
following (as with Proposition~\ref{prop:examplesqaw}, note the occurrence of $\log_0$, i.e.\ of
the identity).

\begin{theorem}\label{thm:detexplicitStieltjes}
Define repeated logarithms $\log_j$ ($j=0,1,2,\ldots$) as in Proposition
\ref{prop:examplesqaw}. For $j=0,1,2,\ldots$ let $a_j>0$ and let $p_j\in\R$ be
such that $p_j=0$ for all sufficiently large $j$. Put
$j_0=\min\{j=0,1,2,\ldots\mid p_j\neq 1\}$ and assume $p_{j_0}<1$.

Let $\{v_1,\ldots,v_n\}$ be a basis of $\Rn$ and let $C=\R_{\geq 0}\cdot v_1+\ldots\R_{\geq
0}\cdot v_n$ be the corresponding positive convex cone. Suppose $\mu$ is a positive Borel
measure on $\Rn$ which is supported in $C$ and such that
\begin{equation*}
\int_{\Vert x\Vert\geq R}\exp\left(\Vert x\Vert^{3/2}\left(\prod_{j=0}^\infty\log_{j}^{p_j}a_j\sqrt{\Vert
x\Vert}\right)^{-1}\right)\,d\mu<\infty,
\end{equation*}
for some $R\geq 0$ which is sufficiently large to ensure that the integrand is
defined. Then $\mu$ is $C$-determinate.
\end{theorem}

As a consequence, if $\mu\in\measures$ is supported in a positive convex cone $C$ as above, and
if one of the functions (tacitly assumed to be equal to $1$ on a sufficiently large compact set)
\begin{align*}
&\exp\left(\frac{\Vert x\Vert^{1/2-\nu}}{a_0}\right)\quad ,\\
&\exp\left(\frac{\sqrt{\Vert x\Vert}}{a_0 \left(\log a_1\Vert
x\Vert\right)^{1+\nu}}\right)\quad ,\\&\exp\left(\frac{\sqrt{\Vert x\Vert}}{a_0 \log
a_1\Vert x\Vert
\left(\log\log a_2\Vert x\Vert\right)^{1+\nu}}\right)\quad ,\\
&\ldots
\end{align*}
has finite integral under $\mu$ for some $a_0,a_1,a_2,\ldots>0$ and $\nu\leq 0$, then $\mu$ is $C$-determinate.

We end this section by showing that in many cases---e.g.\ if the measure is absolutely
continuous with respect to Lebesgue measure---the conclusion of $C$-determinacy in the Theorems~\ref{thm:CarlemanStieltjes}, \ref{thm:maintheoremStieltjes}, \ref{thm:tensorStieltjes},
\ref{thm:radialStieltjes} and \ref{thm:detexplicitStieltjes} can be strengthened to determinacy.
It is sufficient to discuss strengthening Theorem~\ref{thm:CarlemanStieltjes}, since this result
implies the others. We assume that the basis $\{v_1,\ldots,v_n\}$ generating the cone $C$ is the
standard basis; this simplifies the discussion and the general case follows from this by a
linear transformation.

To start with, note that the hypotheses in Theorem~\ref{thm:CarlemanStieltjes} imply that the
marginal distributions of $\mu$ are determinate in the sense of Stieltjes. Recall (see \cite[p.\
481]{Chihara}) that in one dimension a measure, which is supported in $[0,\infty)$ and which is
determinate in the sense of Stieltjes, is actually determinate if its support does not contain
$0$ and/or its support is not equal to a discrete unbounded set. Therefore, if the support of
each marginal distribution satisfies this condition, all marginal distributions are actually
determinate. The results in \cite{Petersen} then imply that $\mu$ itself is determinate. To
summarize:

\emph{Let $\{v_1,\ldots,v_n\}$, $C$ and $\mu$ be as in Theorem~\ref{thm:CarlemanStieltjes}, \ref{thm:maintheoremStieltjes}, \ref{thm:tensorStieltjes},
\ref{thm:radialStieltjes} or \ref{thm:detexplicitStieltjes}. Define marginal distributions of
$\mu$ in terms of the projections corresponding to the basis $\{v_1,\ldots,v_n\}$. If the
support of each of these marginal distributions does not contain $0$ and/or is not equal to a
discrete unbounded set, then $\mu$ is determinate.}

\section{Closing remark}\label{sec:closing}

As mentioned in Section~\ref{sec:qaw} there exist quasi-analytic weights on
$\Rn$ which are quasi-analytic with respect to a unique basis (up to scaling). For $n
\geq 2$ the demonstration of this phenomenon in \cite{deJeu} is
based on the construction of $n$ strictly positive logarithmically convex
sequences $\{M_j(m)\}_{m=1}^\infty$ ($j=1,\ldots,n$) such that
\begin{equation*}
\sum_{m=1}^\infty M_j(m)^{-1/m}=\infty\quad(j=1,\ldots,n),
\end{equation*}
but
\begin{equation*}
\sum_{m=1}^\infty \left(\max (M_{j_1}(m),M_{j_2}(m))\right)^{-1/m}<\infty\quad(1\leq j_1\neq j_2\leq n).
\end{equation*}
Sequences satisfying the first of these equations also figure in Theorem
\ref{thm:DCtheorem}

Now Theorem~\ref{thm:DCtheorem} can evidently be formulated with respect to any basis of $\Rn$,
leading naturally to the notion of quasi-analytic classes with respect to bases. It is an
interesting question whether there then exists an analogue of the aforementioned phenomenon for
quasi-analytic weights. More precisely: can one, perhaps using sequences as above satisfying
\emph{both} equations, establish the existence of smooth functions that are in a quasi-analytic
class with respect to a
\emph{unique} basis, up to scaling? If so, then in view of the proof of the
extended Carleman theorem, additional argumentation could conceivably lead to
the construction of multidimensional measures to which the extended Carleman
theorem applies, but applies with \emph{only one} basis, again up to scaling.
Such measures would then have a distinguished set of marginal distributions.

%%%%%%%%%%%%%%   End of the text   %%%%%%%%%%%%%%%%%%%%%%%%%%%%%%%%%%%%%%%%%%%%%%%%%%%%

%%%%%%%%%%%%%%   Bibliography   %%%%%%%%%%%%%%%%%%%%%%%%%%%%%%%%%%%%%%%%%%%%%%%%%%%%%%%

%%%%%%%%%%%%%%%%%%%%%%%%%%%%%%%%%%%%%%%%%%%%%%%%%%%%%%%%%%%%%%%%%%%%%%%%%%%%%%%%%%%%%%%

%%%%%%%%%%%%%%   End of the document   %%%%%%%%%%%%%%%%%%%%%%%%%%%%%%%%%%%%%%%%%%%%%%%%

\begin{thebibliography}{99}

\bibitem{Berg1}
Berg, C.\ (1995). Recent results about moment problems. {\it Probability measures on groups and
related structures, XI (Oberwolfach, 1994)}, 1--13. World Sci.\ Publishing, River Edge, New
Jersey.

\bibitem{Berg2}
Berg, C.\ (1996). Moment problems and polynomial approximation. 100 ans apr\`es
Th.-J.~Stieltjes. {\it Ann.\ Fac.\ Sci.\ Toulouse Math.\ }{\bf 6}, Special issue, 9--32.

\bibitem{BergChristensen1}
Berg, C.\ and Christensen, J.P.R.\ (1981). Density questions in the classical theory of moments.
{\it Ann.\ Inst.\ Fourier (Grenoble)} {\bf 31} 99--114.

\bibitem{BergChristensen2}
Berg, C.\ and Christensen, J.P.R.\ (1983). Exposants critiques dans le probl\`eme des moments.
{\it C.R.\ Acad.\ Sci.\ Paris S\'er.\ I Math.\ }{\bf 296} 661--663.

\bibitem{Chihara}
Chihara, T.S.\ (1968). On indeterminate Hamburger moment problems. {\it Pacific J.\ Math.\ }{\bf
27} 475--484.

\bibitem{DunfordSchwartz}
Dunford, N.\ and Schwartz, J.T.\ (1957). {\it Linear Operators. Part I}. John Wiley \& Sons,
Inc., New York.

\bibitem{Fuglede}
Fuglede, B.\ (1983). The multidimensional moment problem. {\it Exposition.\ Math.\ }{\bf 1}
47--65.

\bibitem{Hoffmann}
Hoffman-Jorgensen, J.\ (preprint 2002). The Moment Problem.

\bibitem{Hryptun}
Hryptun, V.G.\ (1976). An addition to a theorem of S.~Mandelbrojt. {\it Ukra\"\i n.\ Mat.\
\v Z.\ }{\bf 28} 841--844. English translation: {\it Ukrainian Math.\ J.\ }{\bf 28} 655--658.

\bibitem{deJeu}
De~Jeu, M.F.E.\ (2001). Subspaces with equal closure. To appear (math.CA/0111015).


\bibitem{Koosis}
Koosis, P.\ (1988). {\it The logarithmic integral. I.} Cambridge University Press, Cambridge.

\bibitem{Lin}
Lin, G.D.\ (1997). On the moment problems. {\it Statist.\ Probab.\ Lett.\ }{\bf 35} 85--90.


\bibitem{Nussbaum}
Nussbaum, A.E.\ (1966). Quasi-analytic vectors. {\it Ark.\ Mat.\ }{\bf 6} 179--191.

\bibitem{PakesHungWu}
Pakes, A.G., Hung, W.-L. and Wu, J.-W. (2001). Criteria for the unique determination of
probability distributions by moments. {\it Aust.\ N.\ Z.\ J.\ Stat.\ } {\bf 43} 101--111.


\bibitem{Petersen}
Petersen, L.C.\ (1982). On the relation between the multidimensional moment problem and the
one-dimensional moment problem. {\it Math.\ Scand.\ }{\bf 51} 361--366.

\bibitem{ShohatTamarkin}
Shohat, J.\ and Tamarkin, J.D.\ (1943). {\it The Problem of Moments.} American Mathematical
Society, Providence.

\end{thebibliography}
\end{document}